\documentclass[12 pt]{article}
\usepackage{amsmath}
\usepackage{amssymb}
\usepackage{amsthm}
\usepackage{xypic}

\title{The Jumping Phenomenon of Hodge Numbers}
\begin{document}
\maketitle {$$\bf Xuanming Ye $$} 
 \vskip 0.5cm

\begin{abstract}\ \rm Let $X$ be a compact complex
manifold, consider a small deformation $\phi: \mathcal{X} \rightarrow B$
of $X$, the dimension of the Dolbeault cohomology groups
$H^q(X_t,\Omega_{X_t}^p)$ may vary under this defromation. This
paper will study such phenomenons by studying the obstructions to
deform a class in $H^q(X,\Omega_X^p)$ with the parameter $t$ and
get the
formula for the obstructions.\\
\end{abstract}

\section{Introduction}
Let $X$ be a compact complex manifold and $\phi:\mathcal{X}\rightarrow B$
be a family of complex manifolds such that $\phi^{-1}(0)=X$. Let
$X_{t}=\phi^{-1}(t)$ denote the fibre of $\phi$ above the point
$t\in B$. We denote by $\mathcal{O}_{X}$ and $\Omega^{p}_{X}$ the
sheaves of germs of $X$ of holomorphic functions and $p$-forms
respectively. Recall $h^{p,q}
=dim_{\mathbb{C}}H^{q}(X,\Omega_{X}^{p})$ and $P_{m}=dim
H^{0}(X,(\Omega_{X}^{n})^{\bigotimes m})$ where
$n=dim_{\mathbb{C}}X$. S.Iitaka proposed a problem whether all
$P_{m}$ are deformation invariants [1]. This problem was solved by
Iku Nakamura in his paper [2], and actually he gave us some
examples of small deformations of complex parallelisable manifold
(by a complex parallelisable manifold we mean a compact complex
manifold with the trivial holomorphic tangent bundle) such that
the hodge numbers of the fibre of the family jump in these
deformations. \vskip 0.1cm In this paper, we will study such
phenomenons from the viewpoint of obstruction theory. More
precisely, for a certain small deformation $\mathcal{X}$ of $X$ parametrized by a basis $B$ and a certain
class $[\alpha]$ of the Dolbeaut cohomology group
$H^{q}(X,\Omega_{X}^{p})$, we will try to find out the obstruction
to extending it to an element of the relative Dolbeaut cohomology
group $H^{q}(\mathcal{X},\Omega_{\mathcal{X}/B}^{p})$. We will call those
elements which have non trivial obstruction the obstructed
elements. \vskip 0.1cm

In $\S2$ we will summarize the results of
 Grauert's
Direct Image Theorems and we will try to explain why we need to
consider the obstructed elements. Actually, we will see that these
elements will play an important role when we study the jumping
phenomenon of Hodge numbers. Because we will see that the
existence of the obstructed elements is a necessary and sufficient
condition for the variation of the Hodge diamond. \vskip 0.1cm

In $\S3$ we will get a formula for the obstruction to the
extension we mentioned above.  \vskip 0.1 cm \noindent {\bf
\textbf{Theorem 3.3}} \ \ \ {\it Let $\pi:\mathcal{X}\rightarrow B$ be a
deformation of $\pi^{-1}(0)=X$, where $X$ is a compact complex
manifold. Let $\pi_{n}:X_{n}\rightarrow B_{n}$ be the $n$th order
deformation of $X$. For arbitrary $[\alpha]$ belongs to
$H^q(X,\Omega^p)$, suppose we can extend $[\alpha]$ to order $n-1$
in $H^q(X_{n-1},\Omega^p_{X_{n-1}/B_{n-1}})$. Denote such element
by $[\alpha_{n-1}]$. The obstruction of the extension of
$[\alpha]$ to $n$th order is given by:
$$ o_{n,n-1}(\alpha)=d_{X_{n-1}/B_{n-1}}\circ \kappa_{n} \llcorner (\alpha_{n-1})+
\kappa_{n} \llcorner \circ d_{X_{n-1}/B_{n-1}}(\alpha_{n-1}),$$
where $\kappa_{n}$ is the $n$th order Kodaira-Spencer class and
$d_{X_{n-1}/B_{n-1}}$ is the relative differential operator of the
$n-1$th order deformation.}

\vskip 0.1cm In $\S4$ we will use this formula to study carefully
the example given by Iku Nakamura, i.e. the small deformation of
the
Iwasama manifold and discuss some phenomenons. \vskip 0.1cm 
\indent {\em Acknowledgement.} The research was partially supported by 
China-France-Russian mathematics collaboration grant, No. 
34000-3275100, from Sun Yat-sen University.  The author would also like 
to thank ENS, Paris for its hospitality during the academic years of 
2005--2007. Last but most, the author would like to thank Professor 
Voisin for her patient helps and
valuable suggestions.

\section{Grauert's Direct Image Theorems and Deformation theory}
\indent

In this section, let us first review some general results of
deformation theory. Let $X$ be a compact complex manifold. The
manifold $X$ has an underlying differential structure, but given
this fixed underlying structure there may be many different
complex structures on $X$. In particular, there might be a range
of complex structures on $X$ varying in an analytic manner. This
is the object that we will study.\vskip
 0.3cm\noindent   {\bf
\textbf{Definition 1.0}}\ \ \  A deformation of $X$ consists of a
smooth proper morphism $\phi:\mathcal{X}\rightarrow B$, where $\mathcal{X}$ and
$B$ are connected complex spaces, and an isomorphism $X\cong
\phi^{-1}(0)$, where $0\in B$ is a distinguished point. We call
$\mathcal{X}\rightarrow B$ a family of complex manifolds. \vskip 0.1cm
Although $B$ is not necessarily a manifold, and can be singular,
reducible, or non-reduced, (e.g.
$B=Spec\,\mathbb{C}[\varepsilon]/(\varepsilon^{2})$), since the
problem we are going to research is the phenomenon of the jumping
of the Dolbeaut cohomology, we may assume that $\mathcal{X}$ and $B$ are
complex manifolds.\vskip 0.1cm In order to study the jumping of
the Dolbeaut cohomology, we need the following important theorem
(one of the Grauert's Direct Image Theorems).\vskip 0.1cm
\noindent {\bf \textbf{Theorem 1.1}} \ \ \ {\it Let $X$, $Y$ be complex
spaces, $\pi: X\rightarrow Y$ a proper holomorphic map. Suppose
that $Y$ is Stein, and let $\mathcal{F}$ be a coherent analytic
sheaf on $X$. Let $Y_{0}$ be a relatively compact open set in
$Y$. Then, there is an integer $N>0$ such that the following hold.\\
I. There exists a complex\\
$$ \mathcal{E}^{\cdot}:...\rightarrow \mathcal{E}^{-1} \rightarrow
\mathcal{E}^{0} \rightarrow ... \rightarrow
\mathcal{E}^{N}\rightarrow 0 $$ of finitely generated locally free
$\mathcal{O}_{Y_{0}}$-modules on $Y_{0}$ such that for any Stein
open set $W\subset Y_{0}$, we have\\
$$
H^{q}(\Gamma(W,\mathcal{E}^{\cdot}))\simeq\Gamma(W,R^{q}\pi_{*}(\mathcal{F}))
\simeq H^{q}(\pi^{-1}(W),\mathcal{F}) \qquad \forall q\in
\mathbb{Z}. $$ \\
II. (Base Change Theorem). Assume, in addition, that $\mathcal{F}$
is $\pi$-flat [i.e. $\forall x\in X$, the stalk
$\mathcal{F}_{_{x}}$ is flat over as a module over
$\mathcal{O}_{Y,\pi(x)}$]. Then, there exists a complex \\
$$\mathcal{E}^{\cdot}: 0 \rightarrow \mathcal{E}^{0}\rightarrow
\mathcal{E}^{1} \rightarrow ... \rightarrow \mathcal{E}^{N}
\rightarrow 0 $$ \\
of finitely generated locally free $\mathcal{O}_{Y_{0}}$-sheaves
$\mathcal{E}^{p}$ with the following property:\\
\indent Let $S$ be a Stein space and $f: S\rightarrow Y$ a
holomorphic map. Let $X^{'}=X\times_{Y} S$ and $f^{'}:
X^{'}\rightarrow X$ and $\pi^{'}: X^{'}\rightarrow S$ be the two
projections. Then, if $T$ is an open Stein subset of
$f^{-1}(Y_{0})$, we have, for all $q\in \mathbb{Z\mathbb{}}$\\
$$ H^{q}(\Gamma(T,f^{*}(\mathcal{E}^{\cdot})))\simeq
\Gamma(T,R^{q}\pi^{'}_{*}(\mathcal{F}^{'}))\simeq
H^{q}(\pi^{'-1}(T),\mathcal{F}^{'}) $$ where
$\mathcal{F}^{'}=(f^{'})^{*}(\mathcal{F})$.}\vskip 0.1cm Let
$X$,$Y$ be complex spaces, $\pi:X\rightarrow Y$ a proper map. Let
$\mathcal{F}$ be a $\pi$-flat coherent sheaf on $X$. For $y\in Y$,
denote by $\mathcal{M}_{y}$ the $\mathcal{O}_{Y}$-sheaf of germs
of holomorphic functions "vanishing at y": the stalk of
$\mathcal{M}_{y}$ at y is the maximal ideal of
$\mathcal{O}_{Y,y}$; that at $t\neq y"$ is $\mathcal{O}_{Y,t}$. We
set $\mathcal{F}(y)=$ analytic restriction of $\mathcal{F}$ to
$\pi^{-1}(y)=\mathcal{F}\bigotimes_{\mathcal{O}_{Y}}(\mathcal{O}_{Y}/\mathcal{M}_{y})$.
Since we just need to study the local properties, we may assume,
in view of Theorem 1.0, part II, that there is a complex\\

$$ \mathcal{E}^{\cdot}: \xymatrix { 0 \ar[r] &
\mathcal{O}_{Y}^{P_{0}} \ar[r]^{d^{0}} & \mathcal{O}_{Y}^{P_{1}}
\ar[r]^{d^{1}} & ... \ar[r]^{d^{N-1}} & \mathcal{O}_{Y}^{P_{N}}
\ar[r]^{d^{N}} &
0 \\
}$$ \\

with the base change property in Theorem 1.0, part II. In
particular, if $y\in Y$, we have\\
$$ H^{q}(\pi^{-1}(y),\mathcal{F}(y))\simeq
H^{q}(\mathcal{E}^{\cdot}\otimes(\mathcal{O}_{Y}/\mathcal{M}_{y})).$$\vskip
0.1cm \indent Apply what we discussed above to our case $\phi:\mathcal{X}
\rightarrow B$, we get the following. There is a complex of vector
bundles on the basis $B$, whose cohomology groups at the point
identifies to the cohomology groups of the fiber $X_{b}$ with
values in the considered vector bundle on $\mathcal{X}$, restricted to
$X_{b}$. Therefore, for arbitrary $p$, there exists a complex of
vector bundles $(E^{\cdot},d^{\cdot})$, such that for arbitrary
$t\in B$,
$H^{q}(X_{t},\Omega_{X_{t}}^{p})=H^{q}(E_{t}^{\cdot})=Ker(d^{q})/Im(d^{q-1})$.\vskip
0.1cm \indent Via a local trivialisation of the bundle $E^{i}$,
the differential of the complex $E^{\cdot}$ are represented by
matrices with holomorphic coefficients, and follows from the lower
semicontinuity of the rank of a matrix with variable coefficients
, it is easy to check that the function
$dim_{\mathbb{C}}Ker(d^{q})$ and $-dim_{\mathbb{C}}Im(d^{q})$ are
upper semicontinuous on $B$. Therefore the function
$dim_{\mathbb{C}}H^{q}(E_{t}^{\cdot})$ is also upper
semicontinuous. It seems that either the increasing of
$dim_{\mathbb{C}}Im(d^{q-1})$ or the decreasing of
$dim_{\mathbb{C}}Ker(d^{q})$ will cause the jumping of
$dim_{\mathbb{C}}H^{q}(E_{t}^{\cdot})$, however, because of the
following exact sequence:
$$ 0 \rightarrow Ker(d^{q})_{t} \rightarrow E_{t}^{q} \rightarrow
Im (d^{q})_{t} \rightarrow 0 \indent \forall t,$$ which means the
variation of $-dim_{\mathbb{C}}Im(d^{q})$ is exactly the variation
of $dim_{\mathbb{C}}Ker(d^{q})$, we just need to consider the
variation of $dim_{\mathbb{C}}Ker(d^{q})$ for all $q$.

In order to study the variation of $dim_{\mathbb{C}}Ker(d^{q})$,
we need to consider the following problem. Let $\alpha$ be an
element of $Ker(d^{q})$ at $t=0$, we try to find out the
obstruction to extending it to an element which belongs to
$Ker(d^{q})$ in a neighborhood of $0$.
Such kind of extending can be studied order by order. Let $\mathcal{E}^{q}_{0}$ be the stalk of the assocaited sheaf of $E^{q}$ at $0$. Let $m_{0}$ be the maximal idea of $\mathcal{O}_{B,0}$. For arbitrary positive intergal $n$, since $d^{q}$ can be represented by matrices with holomorphic coefficients, it is not difficult to check $d^{q}(\mathcal{E}^{q}_{0}\otimes_{\mathcal{O}_{B,0}}m_{0}^{n}) \subset \mathcal{E}^{q+1}_{0}\otimes_{\mathcal{O}_{B,0}}m_{0}^{n}$. Therefore the complex of the vector bundles $(E^{\cdot},d^{\cdot})$ induces the following complex:\\
$$ 0 \rightarrow \mathcal{E}^{0}_{0}\otimes_{\mathcal{O}_{B,0}}\mathcal{O}_{B,0}/m_{0}^{n}\stackrel{d^{0}}{\rightarrow} \mathcal{E}^{1}_{0}\otimes_{\mathcal{O}_{B,0}}\mathcal{O}_{B,0}/m_{0}^{n}\stackrel{d^{1}}{\rightarrow}...\stackrel{d^{N-1}}{\rightarrow} \mathcal{E}^{N}_{0}\otimes_{\mathcal{O}_{B,0}}\mathcal{O}_{B,0}/m_{0}^{n}\stackrel{d^{N}}{\rightarrow} 0.$$
\vskip
 0.3cm\noindent   {\bf
 \textbf{Definition 2.2}}\ \ \  Those elements of $H^{\cdot}(E_{0}^{\cdot})$ which can not be extended are called {\em the first class obstructed elements}.\\
\indent Next, we will show the obstructions of the extending we mentioned above. For simplicity, my may assume that $dim_{\mathbb{C}}B=1$, suppose $\alpha$ can be extended to an element $\alpha_{n-1}$ such that $j^{n-1}_{0}(d^{q}(\alpha_{n-1}))(t)=0$, then $\alpha_{n-1}$ can be considered as the $n-1$ order extension of $\alpha$. Here $j^{n-1}_{0}(d^{q}(\alpha_{n-1}))(t)$ is the $n-1$ jet of $d^{q}(\alpha_{n-1})$ at $0$.\\
\indent Define a map  $o^{q}_{n}: H^{q}(\mathcal{E}^{\cdot}_{0}\otimes_{\mathcal{O}_{B,0}}\mathcal{O}_{B,0}/m_{0}^{n}) \rightarrow H^{q+1}(E_{0}^{\cdot})$ by \\
$$ [\alpha_{n-1}] \longmapsto [j^{n}_{0}(d^{q}(\alpha_{n-1}))(t)/t^{n}] .$$
\indent At first, we need to check $o^{q}_{n}$ is well defined. So
we need to show that $[j^{n}_{0}(d^{q}(\alpha_{n-1}))(t)/t^{n}]$
is $d^{q+1}$-closed. Via a local trivialization of the bundles
$E^{i}$, the differentials of the complex $E^{\cdot}$ are
represented by matrices with holomorphic coefficients, and from
the lower semi-continuity of the rank of a matrix with variable
coefficients, we may assume that there always exists
$(\sigma^{q+1}_{1},...,\sigma^{q+1}_{l})$ which are sections of
$E^{q+1}$ such that
$(\sigma^{q+1}_{1}|_{t=0},...,\sigma^{q+1}_{l}|_{t=0})$ form a
basis of $Ker(d^{q+1}:E^{q+1}_{0}\rightarrow E^{q+2}_{0})$ and
$Ker(d^{q+1}:E^{q+1}\rightarrow E^{q+2})\subset
Span\{\sigma^{q+1}_{j}\}$. So we can write
$d^{q}(\alpha_{n-1})=\sum_{j}f_{j}\sigma^{q+1}_{j}$.

Since $j^{n-1}_{0}(d^{q}(\alpha_{n-1}))(t)$=0, we have $f_{j}=0$
and $\frac{\partial^{i}f}{\partial t^{i}}=0$, $i=1..n-1$.
$$\frac{\partial^{n}}{\partial
t^{n}}(d^{q}(\alpha_{n-1}))|_{t=0}=\sum_{j}\frac{\partial^{n}
f_{j}}{\partial
t^{n}}\sigma_{j}^{q+1}|_{t=0}+...+\sum_{j}f_{j}\frac{\partial^{n}}{\partial
t^{n}}(\sigma_{j}^{q+1})|_{t=0}=\sum_{j}\frac{\partial^{n}
f_{j}}{\partial
t^{n}}\sigma_{j}^{q+1}|_{t=0},$$ therefore\\
$$d^{q+1}(\frac{\partial^{n}}{\partial
t^{n}}(d^{q}(\alpha))|_{t=0})=d^{q+1}(\sum_{j}\frac{\partial^{n}
f_{j}}{\partial t^{n}}\sigma_{j}^{q+1}|_{t=0})=0,$$ which means
$\frac{\partial^{n}}{\partial t^{n}}(d^{q}(\alpha_{n-1}))|_{t=0}$
is $d^{q+1}$-closed. \vskip 0.1cm \indent Next we are going to
show that the equivalent class of $\frac{\partial^{n}}{\partial
t^{n}}(d^{q}(\alpha_{n-1}))|_{t=0}$ in $H^{q+1}(E_{0}^{\cdot})$
depends only on $j^{n-1}_{0}(\alpha_{n-1})(t)$. Let
$(\sigma_{1}^{q},...,\sigma_{k}^{q})$ be a bases of $E^{q}$, we
only need to show that if $j^{n-1}_{0}(\alpha_{n-1})(t)=0$, then
$\frac{\partial^{n}}{\partial t^{n}}(d^{q}(\alpha_{n-1}))|_{t=0}$
belongs to $Im(d^{q}:E^{q}\rightarrow E^{q+1})$. Indeed, we can
write $\alpha_{n-1}=\sum_{j}f_{j}\sigma^{q}_{j}$ while
$f_{j}(0)=0, \frac{\partial^{i}f}{\partial t^{i}}=0$, $i=1...n-1$,
then,
$$ \frac{\partial^{n}}{\partial t^{n}}(d^{q}(\alpha_{n-1}))=\frac{\partial^{n}}{\partial
t^{n}}(\sum_{i}f_{i}d^{q}(\sigma^{q}_{i}))=\sum_{i}\frac{\partial^{n}
f_{i}}{\partial
t^{n}}d^{q}(\sigma^{q}_{i})+...+\sum_{i}f_{i}\frac{\partial^{n}}{\partial
t^{n}}(d^{q}(\sigma^{q}_{i})).$$ Therefore,
$\frac{\partial^{n}}{\partial
t^{n}}(d^{q}(\Omega))|_{t=0}=\sum_{i}\frac{\partial^{n}
f_{i}}{\partial t^{n}}d^{q}(\sigma^{q}_{i})|_{t=0}$, which belongs
to $Im(d^{q}:E^{q}\rightarrow E^{q+1})$. \vskip 0.1cm \indent At
last, we are going to show that the equivalent class of
$\frac{\partial^{n}}{\partial t^{n}}(d^{q}(\Omega))|_{t=0}$ in
$H^{q+1}(E_{0}^{\cdot})$ depends only on the equivalent class of
$\alpha_{n-1}$ in
$H^{q}(\mathcal{E}^{\cdot}_{0}\otimes_{\mathcal{O}_{B,0}}\mathcal{O}_{B,0}/m_{0}^{n})$.
Actually, we only need to show that if $\alpha_{n-1}$ belongs to
$Im(d^{q-1}:\mathcal{E}^{q-1}_{0}\otimes_{\mathcal{O}_{B,0}}\mathcal{O}_{B,0}/m_{0}^{n}
 \rightarrow \mathcal{E}^{q}_{0}\otimes_{\mathcal{O}_{B,0}}\mathcal{O}_{B,0}/m_{0}^{n})$, we will have $\frac{\partial^{n}}{\partial
t^{n}}(d^{q}(\alpha_{n-1}))|_{t=0}$ belongs to
$Im(d^{q}:E^{q}\rightarrow E^{q+1})$. In fact, let
$\alpha_{n-1}^{'}=d^{q-1}(\sum_{j}f_{j}\sigma_{j}^{q-1})$ such
that
$j^{n-1}_{0}(\alpha_{n-1}^{'})(t)=j^{n-1}_{0}(\alpha_{n-1})(t)$. From the discussion above, we have\\
$$ \frac{\partial^{n}}{\partial t^{n}}(d^{q}(\alpha_{n-1}))|_{t=0}=
\frac{\partial^{n}}{\partial
t^{n}}(d^{q}(\alpha_{n-1}^{'}))|_{t=0}=\frac{\partial^{n}}{\partial
t^{n}}(d^{q}(d^{q-1}(\sum_{j}f_{j}\sigma_{j}^{q-1})))=0$$ in
$H^{q}(E^{\cdot})$. \vskip 0.1cm \noindent {\bf \textbf{Remark}} \
\ \ It seems that $j^{n}_{0}(d^{q}(\alpha_{n-1}))(t)/t^{n}$
depends on the connection of $E^{q+1}$. But, by using an induction
argument, it is not difficult to prove that if
$j^{i}_{0}(d^{q}(\alpha_{n-1}))(t)=0, \forall i<n$, then
$j^{n}_{0}(d^{q}(\alpha_{n-1}))(t)$ is independent of the choice
of the connection of $E^{q+1}$.\\
\indent There is natural a map 
$\rho^{q}_{i}:  H^{q}(E_{0}^{\cdot}) \rightarrow H^{q}(\mathcal{E}^{\cdot}_{0}\otimes_{\mathcal{O}_{B,0}}\mathcal{O}_{B,0}/m_{0}^{i+1})  $ given by
 $$[\sigma]  \longmapsto  [t^{i}\sigma], \forall [\sigma] \in H^{q}(E_{0}^{\cdot}).$$
Denote the map  $\rho^{q+1}_{i} \circ o^{q}_{n}: H^{q}(\mathcal{E}^{\cdot}_{0}\otimes_{\mathcal{O}_{B,0}}\mathcal{O}_{B,0}/m_{0}^{n}) \rightarrow H^{q+1}(\mathcal{E}^{\cdot}_{0}\otimes_{\mathcal{O}_{B,0}}\mathcal{O}_{B,0}/m_{0}^{i+1}), \forall i \leq n$ by $o^{q}_{n,i}.$ 
\vskip 0.1cm \indent Next we will show that, for arbitrary $i$, $0<i\leq n$,  $\alpha_{n-1}$ can be
extended to $\alpha_{n}$ which is the $n$th order extension of $\alpha$ such that $j^{i-1}_{0}(\alpha_{n}-\alpha_{n-1})(t)=0$ if and only if
$o_{n,n-i}^{q}([\alpha_{n-1}])$ is trivial. For necessarity, $(\alpha_{n}-\alpha_{n-1})(t)/t^{i}$ is supposed to be the preimage of $o_{n,n-i}^{q}([\alpha_{n-1}])$, so $o_{n,n-i}^{q}([\alpha_{n-1}])$ is trivial. Therefore we just need to check
whether it is sufficient. In fact, if $o_{n,n-i}^{q}([\alpha_{n-1}])$
is trivial, then there exists a section $\beta$ of $\mathcal{E}^{\cdot}_{0}\otimes_{\mathcal{O}_{B,0}}\mathcal{O}_{B,0}/m_{0}^{i+1}$ such that
$d^{q}(\beta)=o_{n,n-i}^{q}([\alpha_{n-1}])$. Then it is not difficult
to check that $\alpha_{n-1}-t^{i}\tilde{\beta}$ is an $n$ th order
extension of $\alpha$ that we need, where $\tilde{\beta}$ is an extension of
$\beta$ in the neighborhood of $0$. Therefore we have the
following proposition. \vskip
 0.1cm\noindent   {\bf \textbf{Proposition 2.3}}\ \ \ {\it Let $\alpha_{n-1}$ be an $n-1$ th
order extension of $\alpha$, for arbitrary $i$, $0<i\leq n$,  $\alpha_{n-1}$ can be
extended to $\alpha_{n}$ which is the $n$th order extension of $\alpha$ such that $j^{i-1}_{0}(\alpha_{n}-\alpha_{n-1})(t)=0$ if and only if
$o_{n,n-i}^{q}([\alpha_{n-1}])=0$.}\\
\vskip 0.3cm \indent In the following, we will show that the
obstructions $o^{q}_{n}([\alpha_{n-1}])$ also play an important
role when we consider about the jumping of
$dim_{\mathbb{C}}Im(d^{q})$. Note that $dim_{\mathbb{C}}Im(d^{q})$
jumps if and only if there exist a section $\beta$ of
$dim_{\mathbb{C}}Ker(d^{q+1})$, such that $\beta_{0}$ is not exact
while $\beta_{t}$ is exact for $t\neq 0$.  \vskip
 0.3cm\noindent   {\bf
 \textbf{Definition 2.4}}\ \ \
 Those nontrivial elements of $H^{\cdot}(E_{0}^{\cdot})$ that can always be extended to a section which is only exact at $t \neq 0$ are called {\em the second class obstructed elements}.\\
\indent Note that if $\alpha$ is exact at $t=0$, it can be
extended to an element which is exact at every point. So the
definition above does not depend on the element of a fixed
equivalent class. \vskip
 0.1cm\noindent   {\bf
 \textbf{Proposition 2.5}}\ \ \ {\it Let $[\beta]$ be an nontrivial element of
$H^{q+1}(E_{0}^{\cdot})$. Then $[\beta]$ is a second class
obstructed element if and only if there exist $n\geq0 $ and
$\alpha_{n-1}$ in
$H^{q}(\mathcal{E}^{\cdot}_{0}\otimes_{\mathcal{O}_{B,0}}\mathcal{O}_{B,0}/m_{0}^{n})$
such that
$o^{q}_{n}([\alpha_{n-1}])=[\beta]$.}
\begin{proof}  \indent If
$o^{q}_{n}([\alpha_{n-1}])=[\beta]$, then
$j^{n}_{0}(d^{q}(\alpha_{n-1}))(t)/t^{n}$ is the extension we
need. On the contrary, if $[\beta]$ is a second class obstructed
element. There exist $\tilde{\beta}$ such that
$\tilde{\beta}_{t}$, $t\neq0$ is exact. Then
$(d^{q})^{-1}(\tilde{\beta})$ is a meromorphic section which has a
pole at $t=0$. Let $n$ be the degree of
$(d^{q})^{-1}(\tilde{\beta})$. Then let $\alpha_{n-1} =
t^{n}(d^{q})^{-1}(\tilde{\beta})$. It is easy to check that
$o^{q}_{n}([\alpha_{n-1}])=[\beta]$.
\end{proof}

\noindent   {\bf \textbf{Proposition 2.6}}\ \ \ {\it 
Let $\alpha_{n-1}$ be an element of 
$H^{q}(\mathcal{E}^{\cdot}_{0}\otimes_{\mathcal{O}_{B,0}}\mathcal{O}_{B,0}/m_{0}^{n})$
such that
$o^{q}_{n}([\alpha_{n-1}]) \neq 0 $. Then there exists $n^{'} \leq n$ and $\alpha^{'}$ be an element of 
$H^{q}(\mathcal{E}^{\cdot}_{0}\otimes_{\mathcal{O}_{B,0}}\mathcal{O}_{B,0}/m_{0}^{n^{'}})$, such that $\rho^{q+1}_{n^{'}-1}\circ o^{q}_{n}([\alpha_{n-1}]) = 
o^{q}_{n^{'},n^{'}-1}([\alpha^{'}]) \neq 0 $}.
\begin{proof}  \indent If
$o^{q}_{n,n-1}([\alpha_{n-1}]) \neq 0, $ then $n^{'}=n$ and $\alpha^{'}=\alpha_{n-1}$. Otherwise, there exists $\alpha^{'}_{1}$, such that
$d^{q}(\alpha^{'}_{1})=\rho^{q+1}_{n^-1}\circ o^{q}_{n}([\alpha_{n-1}])$. Note that $o^{q}_{n-1,n-2}([\alpha^{'}_{1}])=\rho^{q+1}_{n-2}\circ o^{q}_{n}([\alpha_{n-1}])=o^{q}_{n,n-2}([\alpha_{n-1}]).$ If we go on step by step as above, we can always get the $n^{'}$ and $\alpha^{'}$ for there is at least one of the 
$o^{q}_{n,i}([\alpha_{n-1}])$ is nontrivial. 
\end{proof}

\indent This proposition tells us that althought $o^{q}_{n}([\alpha_{n-1}]) \neq 0 $ does not mean that $o^{q}_{n,n-1}([\alpha_{n-1}]) \neq 0,$
 we can always find  $\alpha^{'}$ such that $o^{q}_{n}([\alpha_{n-1}])$ comes from obstuctions like $o^{q}_{n,n-1}([\alpha^{'}])$. Therefore we can get the following corollary immediately from Proposition 2.5 and Proposition 2.6.\\
 \noindent {\bf \textbf{Corollary
2.7}} \ \ \ {\it Let $[\beta]$ be an nontrivial element of
$H^{q+1}(E_{0}^{\cdot})$. Then $[\beta]$ is a second class
obstructed element if and only if there exist $n\geq0 $ and
$\alpha_{n-1}$ in
$H^{q}(\mathcal{E}^{\cdot}_{0}\otimes_{\mathcal{O}_{B,0}}\mathcal{O}_{B,0}/m_{0}^{n})$
such that
$o^{q}_{n,n-1}([\alpha_{n-1}])=\rho^{q+1}_{n-1}([\beta])$.}

\indent Let us come back to our problem, suppose $\alpha$ can be extended to an element $\alpha_{n-1}$ such that $j^{n-1}_{0}(d^{q}(\alpha_{n-1}))(t)=0$, since what we care is whether $\alpha$ can be extended to an element which belongs to
$Ker(d^{q})$ in a neighborhood of $0$. So, if we have an $n$th order extension $\alpha_{n}$ of $\alpha$, it is not necessary that
$j^{i-1}_{0}(\alpha_{n}-\alpha_{n-1})(t)=0, \forall i, 1<i<n.$ What we need is just $j^{0}_{0}(\alpha_{n}-\alpha_{n-1})(t)=0$ which means $\alpha_{n}$ is an extension of $\alpha$. So the ``real'' obstructions come from $o_{n,n-1}^{q}([\alpha_{n-1}])$. Since these obstructions is so important when we consider the problem of variation of hodge numbers, we will
try to find out an explicit calculation for such obstructions in next section.

\section{The Formula for the Obstructions}
\indent We are going to prove in this section an explicit formula (Theorem 3.3) for the abstract obstructions described above. Let
$\pi:\mathcal{X}\rightarrow B$ be a deformation of $\pi^{-1}(0)=X$, where
$X$ is a compact complex manifold. For every integer $n\geq 0$,
denote by $B_{n}=Spec\,\mathcal{O}_{B,0}/m_{0}^{n+1}$ the
$n$th order infinitesimal neighborhood of the closed point $0\in
B$ of the base $B$. Let $X_{n}\subset \mathcal{X}$ be the complex space
over $B_{n}$. Let $\pi_{n}:X_{n}\rightarrow B_{n}$ be the $n$th
order deformation of $X$. In order to study the jumping phenomenon
of Dolbeaut cohomology groups, for arbitrary $[\alpha]$ belongs to
$H^{q}(X,\Omega^{p})$, suppose we can extend $[\alpha]$ to order
$n-1$ in $H^q(X_{n-1},\Omega^p_{X_{n-1}/B_{n-1}})$. Denote such
element by $[\alpha_{n-1}]$. In the following, we try to find out
the obstruction of the extension of $[\alpha_{n-1}]$ to $n$th
order. Denote $\pi^{*}(m_{0})$ by $\mathcal{M}_{0}$. Consider the exact sequence\\
$$ 0 \rightarrow \mathcal{M}_{0}^{n} / \mathcal{M}^{n+1}_{0}\otimes \Omega_{X_{0}/B_{0}}^{p}
\rightarrow \Omega_{X_{n}/B_{n}}^{p} \rightarrow
\Omega_{X_{n-1}/B_{n-1}}^{p}\rightarrow 0$$ which induces a long
exact
sequence\\
$$ 0 \rightarrow H^0(X,\mathcal{M}_{0}^{n} / \mathcal{M}^{n+1}_{0}\otimes
\Omega_{X_{0}/B_{0}}^{p})\rightarrow H^0(X_{n},
\Omega_{X_{n}/B_{n}}^{p})\rightarrow
H^0(X_{n-1},\Omega_{X_{n-1}/B_{n-1}}^{p})
$$
$$ \rightarrow H^1(X,\mathcal{M}_{0}^{n} / \mathcal{M}^{n+1}_{0}\otimes
\Omega_{X_{0}/B_{0}}^{p}) \rightarrow ... .$$ The obstruction for
$[\alpha_{n-1}]$ comes from the non trivial image of the
connecting homomorphism
$\delta^{*}:H^q(X_{n-1},\Omega_{X_{n-1}/B_{n-1}}^{p})\rightarrow
H^{q+1}(X,\mathcal{M}_{0}^{n} / \mathcal{M}^{n+1}_{0}\otimes
\Omega_{X_{0}/B_{0}}^{p})$. We will calculate it by
$\breve{C}ech$ calculation. \\
\vskip 0.1cm  Cover $X$ by open sets $U_{i}$ such that, for
arbitrary $i$, $U_{i}$ is small enough. More precisely,
$U_{i}$ is stein and the following exact sequence splits\\
$$ 0 \rightarrow \pi_{n}^{*}(\Omega_{B_{n}})(U_{i})
\rightarrow \Omega_{X_{n}} (U_{i})\rightarrow
\Omega_{X_{n}/B_{n}}(U_{i})\rightarrow 0.$$ So we have a map
$\varphi_{i}:  \Omega_{X_{n}/B_{n}} (U_{i})\rightarrow
\Omega_{X_{n}}(U_{i})$, such that, $\varphi_{i}(
\Omega_{X_{n}/B_{n}} (U_{i}))\oplus
\pi_{n}^{*}(\Omega_{B_{n}})(U_{i}) \cong \Omega_{X_{n}} (U_{i}) $.
Denote by $\iota_{i}$, $\iota^{-1}_{i}$ the inclusion from
$\pi_{n}^{*}(\Omega_{B_{n}})(U_{i})$ to $\Omega_{X_{n}} (U_{i})$
and its inverse. Define $d^{i}_{X_{n}/B_{n}}$ by $ \varphi_{i}
\circ d_{X_{n}/B_{n}}\circ\varphi_{i}^{-1}$ and $d^{i}_{B_{n}}$ by
$ \iota_{i} \circ d_{B_{n}}\circ\iota^{-1}_{i}$. Then it determines a
local decomposition of the exterior differentiation $d_{X_{n}}$ in
$\Omega^{\bullet}_{X_{n}}$\\
$$ d_{X_{n}}=d^{i}_{B_{n}}+d^{i}_{X_{n}/B_{n}}.  $$
Denote the set of alternating $q$-cochains $\beta$ with values in
$\mathcal{F}$ by $\mathcal{C}^{q}(\mathbf{U},\mathcal{F})$, i.e.
to each $q+1$-tuple, $i_{0}<i_{1}...< i_{q}$, $\beta$ assigns a
section $\beta(i_{0},i_{1},..., i_{q})$ of $\mathcal{F}$ over
$U_{i_{0}}\cap U_{i_{1}} \cap ... \cap U_{i_{q}}$. \vskip 0.1cm
Let us still using $\varphi_{i}$ denote the following map,
\begin{eqnarray*}
\varphi_{i}: \pi_{n}^{*}(\Omega^{r}_{B_{n}}) \wedge
\Omega^{p}_{X_{n}/B_{n}}(U_{i}) & \rightarrow &
\Omega^{p+r}_{X_{n}}(U_{i})\\
\varphi_{i}(\omega_{i_{1}}\wedge...\wedge\omega_{i_{r}}\wedge\beta_{j_{1}}\wedge...\wedge\beta_{j_{p}})
 & = &
\omega_{i_{1}}\wedge...\wedge\omega_{i_{r}}\wedge\varphi_{i}(\beta_{j_{1}})\wedge...\wedge\varphi_{i}(\beta_{j_{p}}).
\end{eqnarray*}

Define
$\varphi:\mathcal{C}^{q}(\mathbf{U},\pi_{n}^{*}(\Omega^{r}_{B_{n}})
\wedge \Omega^{p}_{X_{n}/B_{n}})
\rightarrow \mathcal{C}^{q}(\mathbf{U},\Omega^{p+r}_{X_{n}})$ by\\
$$ \varphi(\beta)(i_{0},i_{1},..., i_{q})=\varphi_{i_{0}}(\beta(i_{0},i_{1},...,
i_{q})) \qquad \forall \beta \in
\mathcal{C}^{q}(\mathbf{U},\pi_{n}^{*}(\Omega^{r}_{B_{n}})
\wedge \Omega^{p}_{X_{n}/B_{n}}),$$ where $i_{0}<i_{1}...< i_{q}$.\vskip 0.1cm Define the
total Lie derivative with respect to $B_{n}$\\
$$ L_{B_{n}}: \mathcal{C}^{q}(\mathbf{U},\Omega^{p}_{X_{n}})
\rightarrow \mathcal{C}^{q}(\mathbf{U},\Omega^{p+1}_{X_{n}})$$ by
$$ L_{B_{n}}(\beta)(i_{0},i_{1},...,
i_{q})=d_{B_{n}}^{i}(\beta(i_{0},i_{1},..., i_{q})) \qquad \forall \beta \in \mathcal{C}^{q}(\mathbf{U},\Omega^{p}_{X_{n}}),$$  where
$i_{0}<i_{1}...< i_{q}$.\vskip 0.1cm
 Define, for each $U_{i}$ the
total interior product with respect to $B_{n}$, $I^{i}:
\Omega^{p}_{X_{n}}(U_{i}) \rightarrow \Omega^{p}_{X_{n}}(U_{i}) $
by
$$ I^{i}(\mu dg_{1}\wedge dg_{2}\wedge...\wedge dg_{p})=\mu \sum_{j=1}^{p}dg_{1}\wedge...\wedge dg_{j-1}\wedge d^{i}_{B_{n}}(g_{j})\wedge dg_{j+1}\wedge...\wedge dg_{p}.$$
When $p=0$, we put $I^{i}=0$.\vskip 0.1cm Define $\lambda:
\mathcal{C}^{q}(\mathbf{U},\Omega^{p}_{X_{n}}) \rightarrow
\mathcal{C}^{q+1}(\mathbf{U},\Omega^{p}_{X_{n}}) $ by
$$ (\lambda
\beta)(i_{0},...,i_{q+1})=(I^{i_{0}}-I^{i_{1}})\beta(i_{1},...,i_{q+1}) \qquad \forall \beta \in \mathcal{C}^{q}(\mathbf{U},\Omega^{p}_{X_{n}}).$$

\noindent {\bf \textbf{Lemma 3.0}} \ \ \ {\it
$$ \lambda \circ \varphi \equiv \delta \circ \varphi - \varphi
\circ \delta$$    $\mathrm{mod}.\, 
\pi_{n}^{*}(\Omega^{2}_{B_{n}})\wedge \Omega^{p-1}_{X_{n}}.$}
\begin{proof}
\vskip 0.1cm \noindent Define $J:
C^{q}(\mathbf{U},\Omega^{p}_{X_{n}/B_{n}})\rightarrow
C^{q}(\mathbf(U),\Omega^{p}_{X_{n}})$ by
$$
(J(\beta))(i_{0},...,i_{q+1})=(-1)(\varphi_{i_{0}}-\varphi_{i_{1}})
(\beta(i_{1},...,i_{q+1}), $$ where $i_{0}<i_{1}<...<i_{q+1}$. For
arbitrary $\beta$ belongs to
$C^{q}(\mathbf{U},\Omega^{p}_{X_{n}/B_{n}})$,
\begin{eqnarray*}
(\delta \circ \varphi (\beta)) (i_{0},...,i_{q+1})& = &
\sum_{j=0}^{q+1}(-1)^{j}\varphi(\beta)(i_{0},...,\widehat{i_{j}},...,i_{q+1}) \\
& = &
\varphi_{i_{1}}(\beta)(i_{1},...,i_{q+1})\\
&&
+\sum_{j=1}^{q+1}(-1)^{j}\varphi_{i_{0}}(\beta)(i_{0},...,\widehat{i_{j}},...,i_{q+1}),
\end{eqnarray*}
 while
\begin{eqnarray*}
(\varphi \circ \delta (\beta))(i_{0},...,i_{q+1})& = &
\varphi(\sum_{j=0}^{q+1}(-1)^{j}(\beta)(i_{0},...,\widehat{i_{j}},...,i_{q+1}))
\\ & = &
\sum_{j=0}^{q+1}(-1)^{j}\varphi_{i_{0}}(\beta)(i_{0},...,\widehat{i_{j}},...,i_{q+1})
.
\end{eqnarray*}
 So we have $\delta \circ \varphi-\varphi \circ \delta=J$.\\
\indent Fix $(i_{0},...,i_{q+1})$ and let
$\omega=\beta(i_{1},...,i_{q+1}).$ We must show that
$(I^{i_{0}}-I_{i_{1}})(\varphi_{i_{1}}(\omega))=(-1)(\varphi_{i_{0}}-\varphi_{i_{1}})(w)$
mod $\pi_{n}^{*}(\Omega^{2}_{B_{n}})\wedge \Omega^{p-1}_{X_{n}}$.
By linearity, we may suppose $\varphi_{i_{1}}(\omega)=\mu
dg_{1}\wedge...\wedge dg_{p}$. Then
\begin{eqnarray*}
\varphi_{i_{0}} & = & \mu
d^{i^{0}}_{X_{n}/B_{n}}(g_{1})\wedge...\wedge
d^{i_{0}}_{X_{n}/B_{n}}(g_{p})\\
& = & \mu (dg_{1}-d^{i^{0}}_{X_{n}/B_{n}}(g_{1}))\wedge...\wedge
(dg_{p}-d^{i_{0}}_{X_{n}/B_{n}}(g_{p}))\\
& = & \mu dg_{1}\wedge...\wedge dg_{p}-\sum_{j=1}^{p}\mu
dg_{1}\wedge...dg_{j-1}\wedge d_{B_{n}}^{i_{0}}(g_{j}\wedge
dg_{j+1}\wedge ...\wedge dg_{p}
\end{eqnarray*}
\indent \indent \indent +terms in
 $\pi_{n}^{*}(\Omega^{2}_{B_{n}})\wedge
\Omega^{p-1}_{X_{n}}$.\\
Thus $\varphi_{i_{0}}\equiv \varphi_{i_{1}}(\omega)-I^{i_{0}}\circ
\varphi_{i_{1}} (\omega)$ $\mathrm{mod}.\,
\pi_{n}^{*}(\Omega^{2}_{B_{n}})\wedge \Omega^{p-1}_{X_{n}} $, and
$I^{i_{1}}\circ \varphi_{i_{1}}=0$. which means $\lambda \circ
\varphi \equiv J$ $\mathrm{mod}.\,
\pi_{n}^{*}(\Omega^{2}_{B_{n}})\wedge \Omega^{p-1}_{X_{n}} $.
\end{proof}

Now we are ready to calculate the formula for the obstructions. Let
$\tilde{\alpha}$ be an element of
$\mathcal{C}^{q}(\mathbf{U},\Omega^{p}_{X_{n}/B_{n}})$ such that
its quotient image in
$\mathcal{C}^{q}(\mathbf{U},\Omega^{p}_{X_{n-1}/B_{n-1}})$ is
$\alpha_{n-1}$. Then $\delta^{*}([\alpha_{n-1}])$=
$[\delta(\tilde{\alpha})]$ which is an element of
$H^{q+1}(X,\mathcal{M}_{0}^{n} / \mathcal{M}^{n+1}_{0}\otimes
\Omega_{X_{0}/B_{0}}^{p})\cong \mathrm{m}_{0}^{n} /
\mathrm{m}^{n+1}_{0}\otimes
H^{q+1}(X,\Omega_{X_{0}/B_{0}}^{p})$.\vskip 0.1cm \indent Denote
$r_{X_{n}}$ the restriction to the complex space $X_{n}$. In
order to give the obstructions an explicit calculation, we
need to consider the following map $ \rho:
H^{q}(X,\mathcal{M}^{n}_{0}/\mathcal{M}^{n+1}_{0}\otimes
\Omega^{p}_{X_{0}/B_{0}}) \rightarrow
H^{q}(X_{n-1},\pi_{n-1}^{*}(\Omega_{B_{n}|B_{n-1}})\wedge\Omega^{p}_{X_{n-1}/B_{n-1}}).
$ which is defined by $\rho[\sigma]=[ \varphi^{-1} \circ
r_{X_{n-1}}\circ L_{B_{n}}\circ \varphi (\sigma)].$ \vskip 0.1cm

\noindent {\bf \textbf{Lemma 3.1}} \ \ \ {\it The map: $\rho:
H^{q}(X,\mathcal{M}^{n}_{0}/\mathcal{M}^{n+1}_{0}\otimes
\Omega^{p}_{X_{0}/B_{0}}) \rightarrow \\
H^{q}(X_{n-1},\pi_{n-1}^{*}(\Omega_{B_{n}|B_{n-1}})\wedge\Omega^{p}_{X_{n-1}/B_{n-1}})
$ is well defined.} \vskip 0.1cm \noindent
\begin{proof} At first, we need to show that if $\sigma$ is closed, then
$\varphi^{-1}\circ r_{X_{n-1}}\circ L_{B_{n}}\circ \varphi
(\sigma)$ is closed, which is equivalent to show that $\delta \circ
r_{X_{n-1}}\circ L_{B_{n}}\circ \varphi (\sigma) \equiv 0$
$\mathrm{mod}.\, \pi_{n-1}^{*}(\Omega^{2}_{B_{n}|B_{n-1}})\wedge
\Omega^{p-1}_{X_{n}|X_{n-1}}$.\vskip 0.1cm
 \indent Note that
$d_{X_{n}} \circ \delta=- \delta \circ d_{X_{n}}$. Then
\begin{eqnarray*}
  \delta \circ r_{X_{n-1}}\circ L_{B_{n}}\circ \varphi
(\sigma) & = & r_{X_{n-1}}\circ \delta \circ L_{B_{n}}\circ
\varphi (\sigma) \\ & = & -r_{X_{n-1}} \circ (\delta \circ
d^{\cdot}_{X_{n}/B_{n}}+d^{\cdot}_{X_{n}/B_{n}} \circ \delta +
L_{B_{n}} \circ \delta) \circ \varphi (\sigma).
\end{eqnarray*}
Since \\
$$ L_{B_{n}} \circ \delta \circ \varphi (\sigma) \equiv L_{B_{n}}
\circ (\delta \circ \varphi -\lambda \circ \varphi) (\sigma)
\equiv L_{B_{n}} \circ \varphi \circ \delta(\sigma)=0$$\\
and \\
$$ r_{X_{n-1}}\circ (\delta \circ
d^{\cdot}_{X_{n}/B_{n}}+d^{\cdot}_{X_{n}/B_{n}} \circ \delta)\circ
\varphi(\sigma)=0,$$ we have $\delta \circ r_{X_{n-1}}\circ
L_{B_{n}}\circ \varphi (\sigma) \equiv 0$ $\mathrm{mod}.\,
\pi_{n-1}^{*}(\Omega^{2}_{B_{n}|B_{n-1}})\wedge
\Omega^{p-1}_{X_{n}|X_{n-1}}$.\vskip 0.1cm \indent Next we need to
show that if $\sigma$ is belongs to
$C^{q}(\mathbf{U},\mathcal{M}^{n}_{0}/\mathcal{M}^{n+1}_{0}\otimes
\Omega^{p}_{X_{0}/B_{0}})$, then $\varphi^{-1}\circ
r_{X_{n-1}}\circ L_{B_{n}}\circ \varphi \circ \delta(\sigma)$ is
exact. In fact, as the calculation above:\\
$r_{X_{n-1}}\circ L_{B_{n}}\circ \varphi \circ \delta(\sigma)
\equiv -r_{X_{n-1}} \circ (\delta \circ
d^{\cdot}_{X_{n}/B_{n}}+d^{\cdot}_{X_{n}/B_{n}} \circ \delta +
\delta \circ L_{B_{n}}) \circ \varphi (\sigma)=-\delta \circ
r_{X_{n-1}} \circ L_{B_{n}} \circ \varphi (\sigma) . $\vskip 0.1cm
\indent 

\end{proof}
\indent In general, the map $\rho$ is not injective. However, as we mentioned at the end of the previous section. The ``real'' obstructions are
$o_{n,n-1}^{q}([\alpha_{n-1}])$, but not $o_{n}^{q}([\alpha_{n-1}])$. So we don't need $\rho$ to be injective. In the following, we will explain
that $\rho([\delta(\tilde{\alpha})])$ is exactly the ``real'' obstructions we need. In fact, $$H^{q}(X_{n-1},\pi_{n-1}^{*}(\Omega_{B_{n}|B_{n-1}})\wedge\Omega^{p}_{X_{n-1}/B_{n-1}})=(\Omega_{B_{n}|B_{n-1}})\otimes_{\mathcal{O}_{B_{n-1}}} H^{q}(X_{n-1},\Omega^{p}_{X_{n-1}/B_{n-1}}).$$
Let $m=dim_{\mathbb{C}}B$, let $t_{i}, i=0...m$ be the local coordinates of $B$. Then $\rho([\delta(\tilde{\alpha})])$ can be written as:
$ \sum_{i=0}^{m} dt_{i} \otimes \tilde{\alpha}_{i}$, where $\tilde{\alpha}_{i} \in H^{q}(X_{n-1},\Omega^{p}_{X_{n-1}/B_{n-1}}).$ For a certain direction $\frac{\partial}{\partial t_{i}},$ suppose $\tilde{\alpha}_{i} \neq 0$. Then by a simple calculation, it is not difficult to check that $\tilde{\alpha}_{i}= constant [\delta(\tilde{\alpha})/t_{i}]$ in $H^{q}(X_{n-1},\Omega^{p}_{X_{n-1}/B_{n-1}}).$ While $[\delta(\tilde{\alpha})/t_{i}]$ is exactly the obstruction $o_{n,n-1}^{q}([\alpha_{n-1}])$ in the direction of $\frac{\partial}{\partial t_{i}}$ we mentioned in the previous section. \vskip 0.1cm
 Now consider the
following exact sequence. The connecting homomorphism  of the
associated long exact sequence gives the Kodaira-Spencer
 class of order $n$ [4 1.3.2],\\
 $$ 0\rightarrow \pi_{n-1}^{*}(\Omega_{B_{n}|B_{n-1}})\rightarrow
 \Omega_{X_{n}|X_{n-1}}\rightarrow
 \Omega_{X_{n-1}/B_{n-1}}\rightarrow 0 .$$
By wedge the above exact sequence with
$\Omega^{p-1}_{X_{n-1}/B_{n-1}}$, we get a new exact sequence. The
connecting homomorphism of such exact sequence gives us a map from
$H^{q}(X_{n-1},\Omega^{p}_{X_{n-1}/B_{n-1}})$ to
$H^{q+1}(X_{n-1},\pi^{*}(\Omega_{B_{n}|B_{n-1}})\wedge\Omega^{p-1}_{X_{n-1}/B_{n-1}})$.
Denote such map by $\kappa_{n}\llcorner$, for such map is
simply the inner product with the Kodaira-Spencer
 class of order $n$. By the definition and simply calculation it is not difficult to proof
the following lemma.\\
 \noindent {\bf \textbf{Lemma 3.2}} \ \ \
 {\it Let $\theta$ be an element of
 $H^{q}(X_{n-1},\Omega^{p}_{X_{n-1}/B_{n-1}})$, let
 $\tilde{\theta}$ be an element of
 $\mathcal{C}^{q}(\mathbf{U},\Omega^{p}_{X_{n}/B_{n}})$ such that
 its quotient image is $\theta$. Then $[\kappa_{n}\llcorner \theta]$
 is equal to $[\varphi^{-1} \circ r_{X_{n-1}} \circ \delta \circ
 \varphi (\tilde{\theta})].$} \vskip 0.1cm

 Let us come back to the problem we discussed, we have\\
 \begin{eqnarray*}
r_{X_{n-1}}\circ L_{B_{n}} \circ \varphi \circ
\delta(\tilde{\alpha}) & \equiv & r_{X_{n-1}} \circ L_{B_{n}}
\circ (\delta \circ \varphi - \lambda \circ
\varphi)(\tilde{\alpha})  \\
& \equiv & r_{X_{n-1}} \circ L_{B_{n}} \circ \delta \circ \varphi
(\tilde{\alpha})  \\
& \equiv & -r_{X_{n-1}} \circ (d^{\cdot}_{X_{n}/B_{n}} \circ
\delta + \delta \circ d^{\cdot}_{X_{n}/B_{n}}+ \delta \circ
L_{B_{n}} ) \circ
\varphi (\tilde{\alpha}) \\
& \equiv & -r_{X_{n-1}} \circ (d^{\cdot}_{X_{n}/B_{n}} \circ
\delta + \delta \circ d^{\cdot}_{X_{n}/B_{n}}) \circ \varphi
(\tilde{\alpha}) \\
&& -  \delta \circ r_{X_{n-1}}\circ L_{B_{n}}) \circ
\varphi (\tilde{\alpha}). \\
\end{eqnarray*}

\noindent Therefore
\begin{eqnarray*}
[r_{X_{n-1}}\circ L_{B_{n}} \circ \varphi \circ
\delta(\tilde{\alpha})] & = & [-r_{X_{n-1}} \circ
(d^{\cdot}_{X_{n}/B_{n}} \circ \delta + \delta \circ
d^{\cdot}_{X_{n}/B_{n}}) \circ \varphi (\tilde{\alpha}) ]\\
& = & -[d^{\cdot}_{X_{n-1}/B_{n-1}} \circ r_{X_{n-1}}\delta \circ
\varphi (\tilde{\alpha}) + r_{X_{n-1}} \circ \delta \circ
d^{\cdot}_{X_{n}/B_{n}} \circ \varphi (\tilde{\alpha}) ]\\
& = &  -[d^{\cdot}_{X_{n-1}/B_{n-1}} \circ \varphi \circ \varphi
^{-1} \circ r_{X_{n-1}}\delta \circ \varphi (\tilde{\alpha})\\
&& + r_{X_{n-1}} \circ \delta \circ \varphi \circ
d_{X_{n}/B_{n}}(\tilde{\alpha}) ] \\
& = & -[\varphi \circ d_{X_{n-1}/B_{n-1}} \circ \varphi ^{-1}
\circ r_{X_{n-1}}\delta \circ \varphi (\tilde{\alpha}) \\
& & +r_{X_{n-1}} \circ \delta \circ \varphi \circ
(\widetilde{d_{X_{n-1}/B_{n-1}}(\alpha_{n-1})}) ]\\
& = & -[d_{X_{n-1}/B_{n-1}} \circ \kappa_{n} \llcorner
\alpha_{n-1}+ \kappa_{n} \llcorner \circ
d_{X_{n-1}/B_{n-1}}(\alpha_{n-1})] .\\
\end{eqnarray*}
 \vskip 0.1cm
\indent From the discussion above, we get the main theorem of this
paper.

 \vskip 0.1 cm \noindent {\bf \textbf{Theorem
3.3}} \ \ \ {\it Let $\pi:\mathcal{X}\rightarrow B$ be a deformation of
$\pi^{-1}(0)=X$, where $X$ is a compact complex manifold. Let
$\pi_{n}:X_{n}\rightarrow B_{n}$ be the $n$th order deformation of
$X$. For arbitrary $[\alpha]$ belongs to $H^q(X,\Omega^p)$,
suppose we can extend $[\alpha]$ to order $n-1$ in
$H^q(X_{n-1},\Omega^p_{X_{n-1}/B_{n-1}})$. Denote such element by
$[\alpha_{n-1}]$. The obstruction of the extension of
$[\alpha]$ to $n$th order is given by:
$$ o_{n,n-1}(\alpha_{n-1})=d_{X_{n-1}/B_{n-1}} \circ \kappa_{n} \llcorner (\alpha_{n-1})+
\kappa_{n} \llcorner \circ d_{X_{n-1}/B_{n-1}}(\alpha_{n-1}),$$
where $\kappa_{n}$ is the $n$th order Kodaira-Spencer class and
$d_{X_{n-1}/B_{n-1}}$ is the relative differential operator of the
$n-1$th order deformation.}\vskip 0.1cm
\indent From the theorem, we can get the following corollary immediately.\\
 \noindent {\bf \textbf{Corollary
3.4}} \ \ \ {\it Let $\pi:\mathcal{X}\rightarrow B$ be a deformation of
$\pi^{-1}(0)=X$, where $X$ is a compact complex manifold. Suppose that
up to order $n$, the $d_{1}$ of the Fr\"{o}licher spectral sequence vanishes. For arbitrary $[\alpha]$ belongs to $H^q(X,\Omega^p)$, it can be 
extended to order $n+1$ in
$H^q(X_{n+1},\Omega^p_{X_{n+1}/B_{n+1}})$}.

\section{An Example}
\indent \indent In this section, we will use the formula in
previous section to study the jumping of the Hodge numbers
$h^{p,q}$ of small deformations of Iwasawa manifold. It was
Kodaira who first calculated small deformations of Iwasawa
manifold [2]. In the first part of this section, let us recall his
result. \vskip 0.1cm \indent Set
\begin{displaymath}
G=\left\{\left( \begin{array}{ccc} 1 & z_{2} & z_{3}\\
                            0 & 1 & z_{1}\\
                            0 & 0 & 1\\
                            \end{array} \right); z_{i}\in
                            \mathbb{C}\right\}\cong
                            \mathbb{C}^{3}
\end{displaymath}\\
\begin{displaymath}
\Gamma=\left\{\left( \begin{array}{ccc} 1 & \omega_{2} & \omega_{3}\\
                            0 & 1 & \omega_{1}\\
                            0 & 0 & 1\\
                            \end{array} \right); \omega_{i}\in
                            \mathbb{Z}+\mathbb{Z}\sqrt{-1}
                            \right\}\\.
\end{displaymath}\\
The multiplication is defined by
\begin{displaymath}
\left( \begin{array}{ccc} 1 & z_{2} & z_{3}\\
                            0 & 1 & z_{1}\\
                            0 & 0 & 1\\
                            \end{array} \right)
                            \left( \begin{array}{ccc} 1 & \omega_{2} & \omega_{3}\\
                            0 & 1 & \omega_{1}\\
                            0 & 0 & 1\\
                            \end{array} \right)=
                            \left( \begin{array}{ccc} 1 & z_{2}+\omega_{2} & z_{3}+\omega_{2}z_{1}+\omega_{3}\\
                            0 & 1 & z_{1}+\omega_{1}\\
                            0 & 0 & 1\\
                            \end{array} \right)
\end{displaymath}\\
$X=G/\Gamma$ is called Iwasawa manifold. We may consider
$X=\mathbb{C}^{3}/\Gamma$. $g\in \Gamma$ operates on
$\mathbb{C}^3$ as follows: \\
$$ z'_{1}=z_{1}+\omega_{1},   \qquad    z'_{2}=z_{2}+\omega_{2},  \qquad     z'_{3}=z_{3}+\omega_{1}z_{2}+\omega_{3}
$$\\
where $g=(\omega_{1},\omega_{2},\omega_{3})$ and $z'=z\cdot g$.
There exist holomorphic $1$-froms
$\varphi_{1},\varphi_{2},\varphi_{3}$ which are linearly
independent at every point on $X$ and are given by \\
$$\varphi_{1}=dz_{1},\qquad \varphi_{2}=dz_{2}, \qquad
\varphi_{3}=dz_{3}-z_{1}dz_{2},$$\\
so that
$$ d\varphi_{1}=d\varphi_{2}=0, \qquad
d\varphi_{3}=-\varphi_{1}\wedge\varphi_{2}.$$\\
On the other hand we have holomorphic vector fields
$\theta_{1},\theta_{2},\theta_{3}$ on $X$ given by
$$ \theta_{1}=\frac{\partial}{\partial z_{1}},\qquad \theta_{2}=
\frac{\partial}{\partial z_{2}}+z_{1}\frac{\partial}{\partial
z_{3}},\qquad \theta_{3}=\frac{\partial}{\partial z_{3}},$$ \\
It is easily seen that \\
$$ [\theta_{1},\theta_{2}]=-[\theta_{2},\theta_{1}]=\theta_{3}, \qquad
[\theta_{1},\theta_{3}]=[\theta_{2},\theta_{3}]=0.$$ in view of
Theorem $3$ in [2], $H^{1}(X,\mathcal{O})$ is spanned by
$\overline{\varphi}_{1},\overline{\varphi}_{2}$. Since $\Theta$ is
isomorphic to $\mathcal{O}^{3}$, $H^{1}(X,TX)$ is spanned by
$\theta_{i}\overline{\varphi}_{\lambda}, i=1,2,3, \lambda=1,2$.\\
\indent The small deformation o f $X$ is given by
$$ \psi(t)=\sum^{3}_{i=1}\sum_{\lambda=1}^{2}t_{i\lambda}\theta_{i}\overline{\varphi}_{\lambda}t-
(t_{11}t_{22}-t_{21}t_{12})\theta_{3}\overline{\varphi}_{3}t^{2}.$$ \\
We summarize the numerical characters of deformations. The
deformations are divided into the following three classes:\vskip
0.1 cm
 \noindent i)  $t_{11}=t_{12}=t_{21}=t_{22}=0$, $X_{t}$ is a parallelisable
manifold. \\
 ii)  $t_{11}t_{22}-t_{21}t_{12}=0$ and
$(t_{11},t_{12},t_{21},t_{22})\neq (0,0,0,0)$, $X_{t}$ is not
parallelisable.\\
 iii)  $t_{11}t_{22}-t_{21}t_{12}\neq 0$, $X_{t}$ is not
parallelisable.\vskip 0.1 cm
\begin{tabular}{|c|c|c|c|c|c|c|c|c|c|}
\hline
 \rule{0pt}{3ex} & $h^{1,0}$ & $h^{0,1}$ & $h^{2,0}$ & $h^{1,1}$ & $h^{0,2}$ &
$h^{3,0}$
& $h^{2,1}$ & $h^{1,2}$ & $h^{3,0}$ \\
\hline i) & 3 & 2 & 3 & 6 & 2 & 1 & 6 & 6 & 1 \\
\hline ii) & 2 & 2 & 2 & 5 & 2 & 1 & 5 & 5 & 1 \\
\hline iii) & 2 & 2 & 1 & 5 & 2 & 1 & 4 & 4 & 1 \\
\hline
\end{tabular} \vskip 0.1 cm
\indent Now let us explain the jumping phenomenon of the Hodge
number by using the obstruction formula. From Corollary 4.3 in
[6], it follows that the Dolbeault cohomology groups are:\vskip
0.1 cm \noindent
\begin{eqnarray*}
H^{0}(X,\Omega)& =
&Span\{[\varphi_{1}],[\varphi_{2}],[\varphi_{3}]\},\\ 
H^{1}(X,\mathcal{O})& = &
Span\{[\overline{\varphi}_{1}],[\overline{\varphi}_{2}]\},
\\
H^{0}(X,\Omega^{2})&= &Span\{[\varphi_{1}\wedge
\varphi_{2}],[\varphi_{2}\wedge \varphi_{3}],[\varphi_{3}\wedge
\varphi_{1}]\}, \\
H^{1}(X,\Omega)&=&Span\{[\varphi_{i}\wedge\overline{\varphi}_{\lambda}]\},i=1,2,3,\lambda=1,2,
\\
 H^{2}(X,\mathcal{O})&=&Span\{[\overline{\varphi}_{2}\wedge
\overline{\varphi}_{3}],[\overline{\varphi}_{3}\wedge
\overline{\varphi}_{1}]\}, \\
H^{0}(X,\Omega^{3})&=&Span\{[\varphi_{1}\wedge \varphi_{2}\wedge
\varphi_{3}]\},\\
 H^{1}(X,\Omega^{2})&=&Span\{[\varphi_{i}\wedge
\varphi_{j}\wedge
\overline{\varphi}_{\lambda}]\},i,j=1,2,3,i<j,\lambda=1,2, \\
 H^{2}(X,\Omega^{1})&=&Span\{[\varphi_{i}\wedge
\overline{\varphi}_{2}\wedge
\overline{\varphi}_{3}],[\varphi_{j}\wedge
\overline{\varphi}_{1}\wedge
\overline{\varphi}_{3}]\},i,j=1,2,3, \\
 H^{3}(X,\mathcal{O})&=&Span\{[\overline{\varphi}_{1}\wedge
\overline{\varphi}_{2}\wedge
\overline{\varphi}_{3}]\},\end{eqnarray*}
 For example, let us first consider $h^{2,0}$, in the ii) class of
 deformation. The Kodaira-Spencer class of the this deformation is $\psi_{1}(t)=\sum^{3}_{i=1}\sum_{\lambda=1}^{2}t_{i\lambda}\theta_{i}\overline{\varphi}_{\lambda}$,
 with $t_{11}t_{22}-t_{21}t_{12}=0$. It is easy to check that
 $o_{1}(\varphi_{1}\wedge \varphi_{2})=\partial(int(\psi_{1}(t))(\varphi_{1}\wedge \varphi_{2})-int(\psi_{1}(t))(\partial(\varphi_{1}\wedge
 \varphi_{2}))=0$, $o_{1}(t_{11}\varphi_{2}\wedge
 \varphi_{3}-t_{21}\varphi_{1}\wedge \varphi_{3})=\partial((t_{11}t_{22}-t_{21}t_{12})\varphi_{3}\wedge
 \overline{\varphi}_{2})=0$, and $o_{1}(\varphi_{2}\wedge
 \varphi_{3})=-t_{21}\varphi_{1}\wedge\varphi_{2}\wedge\overline{\varphi}_{1}-
 t_{22}\varphi_{1}\wedge\varphi_{2}\wedge\overline{\varphi}_{2}$, $o_{1}(\varphi_{1}\wedge
 \varphi_{3})=-t_{11}\varphi_{1}\wedge\varphi_{2}\wedge\overline{\varphi}_{1}-
 t_{21}\varphi_{1}\wedge\varphi_{2}\wedge\overline{\varphi}_{2}$.
 Therefore, we have shown that for an element of the subspace
 $Span\{[\varphi_{1}\wedge \varphi_{2}],[t_{11}\varphi_{2}\wedge
 \varphi_{3}-t_{21}\varphi_{1}\wedge \varphi_{3}]\}$, the first
 order obstruction is trivial, while, since
 $(t_{11},t_{12},t_{21},t_{22})\neq (0,0,0,0)$, at least one of
 the obstruction $o_{1}(\varphi_{2}\wedge
 \varphi_{3})$, $o_{1}(\varphi_{1}\wedge
 \varphi_{3})$ is non trivial which partly explain why the Hodge
 number $h^{2,0}$ jumps from 3 to 2. For another example, let us
 consider $h^{1,2}$, in the ii) class of deformation. It is easy
 to check that for an element of the subspace (the dimension of
 such a subspace is 5)
 $Span\{[\varphi_{i}\wedge\overline{\varphi}_{\lambda}\wedge\overline{\varphi}_{3}],
 [t_{12}\varphi_{3}\wedge\overline{\varphi}_{2}\wedge\overline{\varphi}_{3}-
 t_{11}\varphi_{3}\wedge\overline{\varphi}_{1}\wedge\overline{\varphi}_{3}]\},i=1,2,\lambda=1,2,$
 the first order obstruction is trivial, while at least one of the
 obstruction
 $o_{1}(\varphi_{3}\wedge\overline{\varphi}_{2}\wedge\overline{\varphi}_{3})$,
 $o_{1}(\varphi_{3}\wedge\overline{\varphi}_{1}\wedge\overline{\varphi}_{3})$
 is non trivial.\vskip 0.1 cm
 \noindent {\bf \textbf{Remark 1}} \ \ \
 It is easy to see that, in the ii) or iii) class of deformation,
 the first order obstruction for any element in $H^{1}(X,\Omega)$
 is trivial.
 The reason of Hodge number $h^{1,1}$'s jumping from 6 to
 5 comes from the existence of the second class obstructed
 elements $o_{1}(\varphi_{3})$. After simple calculation, it is not difficult to get
 the structure equation of $X_{t}, t\neq 0$.

\begin{displaymath}
\left\{ \begin{array}{ll}
d\varphi_{1}=0, \\ d\varphi_{2}=0, \\
d\varphi_{3}=-\varphi_{1}\wedge \varphi_{2}+ t o_{1}(\varphi_{3}),
\qquad i=1,2,\lambda=1,2,
\end{array} \right.
\end{displaymath}
which can be considered an example of proposition 2.5.\\
\noindent {\bf \textbf{Remark 2}} \ \ \ From the example we
discussed above, it is not difficult to find out the following
fact. Let $X$ be an non-K$\ddot{a}$hler nilpotent complex
parallelisable manifold whose dimension is more than 2, and $\phi:
\mathcal{X}\rightarrow B$ be the versal deformation family of $X$. Then
the Hodge number $h^{1,0}$ will jump in a neighborhood of $0\in
B$. In fact, let $\varphi_{i}, i=1...n, n=dim_{\mathbb{C}}(X)$ be
the linearly independent holomorphic 1-forms of $X$. By the
theorem 3 of [2], $H^{1}(X,\mathcal{O})$ is spanned by a subset of
$\{\overline{\varphi}_{i}\}, i=1..n$. So we have
$\partial:H^{1}(X,\mathcal{O})\rightarrow H^{1}(X,\Omega)$ is
trivial, which means one term of the first order obstruction of
the holomorphic 1-forms vanishes. Let $\theta_{i}, i=1...n$ be the
dual of $\varphi_{i}$, which are linearly independent holomorphic
vector fields. Since $X$ is non-K$\ddot{a}$hler, which means $X$
is not a torus, there exists $\varphi_{i}$ such that
$\partial\varphi_{i}\neq 0$. Since $X$ is nilpotent, there exist
$\varphi_{j}$ such that $\partial\varphi_{j}= 0$. Assume that
$\partial\varphi_{i}=A\varphi_{k}\wedge\varphi_{l}+...$ with
$A\neq 0$. Consider $\theta_{k}\overline{\varphi}_{j}$ in
$H^{1}(X,TX)$. It is easy to check that
$o_{1}(\partial\varphi_{i},\theta_{k}\overline{\varphi}_{j})\neq
0$.


\begin{thebibliography}{[6]}
\rm
\bibitem{} S. Iitaka, Plurigenera and classification of algebraic varieties, Sugaku 24
           (1972), 14-27.
\bibitem{} Nakamura, I(1975). Complex parallelisable manifolds and their small deformations, J.Differential Geom. 10,
           85-112.
\bibitem{} C. Voisin, {\it Hodge Theory and Complex Algebraic Geometry I}, Cambridge University Press 2002.
\bibitem{} C. Voisin, {\it Sym\'{e}trie miroir}, Soci\'{e}t\'{e} Math\'{e}matique de France, Paris, 1996. 
\bibitem{} Bell S. and Narasimhan R., Proper holomorphic mappings of complex spaces,
           Encyclopedia of Mathematical Sciences, Several Complex Variables VI,
           Springer Verlag, pp. 1-38, 1991.
\bibitem{} Cordero, L. A., Fern\'{a}ndez, Gray, A. and Ugate, L.(1999). Fr\"{o}licher Spectral Sequence of Compact Nilmanifolds with Nilpotent     Complex Structure. New developments in differential geometry, Budapest 1996, 77-102, Kluwer Acad. Publ., Dordrecht.
\end{thebibliography}
\end{document}